\documentclass[reqno]{amsart}
\usepackage{amsthm,amsmath,amssymb,amscd,xy}
\usepackage{textcomp}
\usepackage{titlesec}
\usepackage[colorlinks,linkcolor=red,anchorcolor=blue,citecolor=green]{hyperref}
\usepackage{anysize}
\usepackage{setspace}
\usepackage{wasysym}
\usepackage{enumitem}
\usepackage{graphicx,color,transparent}
\usepackage{caption}
\usepackage{subcaption}

\newtheorem{theorem}{Theorem}
\newtheorem{lemma}[theorem]{Lemma}
\newtheorem{corollary}[theorem]{Corollary}
\newtheorem*{theorem*}{Theorem}
\newtheorem*{proposition*}{Proposition}

\theoremstyle{definition}

\theoremstyle{remark}

\newcommand\numberthis{\addtocounter{equation}{1}\tag{\theequation}}

\newcommand{\PSL}{\rm{PSL}}

\newcommand{\tr}{\rm{tr}}

\newcommand{\SL}{\rm{SL}}

\newcommand{\Z}{\mathbb{Z}}
\newcommand{\C}{\mathbb{C}}

\allowdisplaybreaks[1]
\titleformat{\section}{\large\bfseries}{\thesection}{1em}{}
\titleformat{\author}{\large\bfseries}{\theauthor}{1em}{}
\marginsize{3cm}{3cm}{3cm}{3cm}
\begin{document}
\title{New identities for small hyperbolic surfaces}
\author{Hengnan Hu and Ser Peow Tan}

\address{Department of Mathematics \\ National University of Singapore \\
 Singapore 119076} \email{mathh@nus.edu.sg; mattansp@nus.edu.sg}

\thanks{
The second author is partially supported by the National University
of Singapore academic research grant R-146-000-186-112. }

\date{}
\maketitle

\begin{abstract}
Luo and Tan gave a new identity for hyperbolic surfaces with/without geodesic boundary in terms of dilogarithms of the lengths of simple closed geodesics on embedded three-holed spheres or one-holed tori in \cite{LuoTan11}. However, the identity was trivial for a hyperbolic one-holed torus with geodesic boundary. In this paper we adapt the argument from Luo and Tan to give an identity for hyperbolic tori with one geodesic boundary or cusp in terms of dilogarithm functions on the set of lengths of simple closed geodesics on the torus. As a corollary, we are also able to express the Luo-Tan identity as a sum over all immersed three-holed spheres $P$ which are embeddings when restricted to the interior of $P$.
\end{abstract}

\section{Introduction}
In the last couple of decades, several interesting identities had been derived by various authors for hyperbolic manifolds and surfaces, notably, those of Basmajian \cite{Basmajian93}, McShane \cite{McShane98,McShane04,McShane06}, Bridgeman \cite{Bridgeman11}, Bridgeman and Kahn \cite{BridgemanKahn10} and Luo and Tan \cite{LuoTan11}. For a survey of the identities, generalizations and their connections, see Bridgeman and Tan \cite{BridgemanTan13b} and also \cite{ BridgemanTan13a}.  For hyperbolic surfaces, the identities derived were for surfaces with geodesic boundary and/or cusps, except for the Luo-Tan identity which works for surfaces with or without boundary. However, the identity derived by Luo-Tan was trivial in the case of hyperbolic one-holed tori as in this case, there was only one embedded simple surface, namely the surface itself. This contrasts for example with McShane's identity which was first proved for hyperbolic tori with one cusp by McShane in his thesis \cite{McShane91} and involves a non-trivial sum over the lengths of all simple  closed geodesics, and subsequently generalized by Mirzakhani for hyperbolic tori $T$ with one geodesic boundary (more generally hyperbolic surfaces with geodesic boundary) in \cite{Mirzakhani07}, and Tan, Wong and Zhang for hyperbolic tori with one cone singularity (more generally hyperbolic surfaces with cone singularities of cone angle $\le \pi$) in \cite{TanWongZhang06}, and more generally for $\SL (2, \C)$ characters of  $\pi_1(T)$ satisfying Bowditch's Q-conditions in \cite{TanWongZhang08, TanWongZhang06b}.

The main aim of this paper is to derive a non-trivial version of the Luo-Tan identity for hyperbolic tori with one geodesic boundary or cusp, and as a corollary, to express the Luo-Tan identity for closed hyperbolic surfaces  as a sum over all immersed geometric three-holed spheres $P$ in the surface which are embeddings when restricted to the interior of $P$. The basic idea is to exploit the close connection between four-holed spheres and one-holed tori, and to use the non-trivial Luo-Tan identity for hyperbolic four-holed spheres with four geodesic boundaries of equal lengths to obtain a corresponding identity for the one-holed tori.  An earlier version of the results here form part of the thesis of the first author \cite{Hu13}. 

Let $T$ be a hyperbolic torus with one geodesic boundary or cusp. Let $$SG_T=\{\gamma ~:~ \gamma ~~\hbox{is a non-peripheral simple closed geodesic on $T$}\}$$  and let $SS_T$ be the set of lengths (counted with multiplicity) of the elements of $SG_T$, we call $SS_T$ the {\it simple length spectrum} of $T$.
The main result of this paper is the following identity for $T$:

\begin{theorem}
\label{thm:oneholed}
For a hyperbolic one-holed torus $T$ with boundary geodesic $K$ of length $k>0$, we have

\begin{equation}
\label{eq:main}
\sum_{b}\Bigg\{\mathcal{L}\left(\frac{\cosh(\frac{k}{2})+1}{\cosh(\frac{k}{2})+\cosh(b)}\right)+ 2 \mathcal{L}\left(\frac{\cosh(\frac{k}{4}+\frac{b}{2})}{\cosh(\frac{k}{4}) e^{\frac{b}{2}}}\right)
-2 \mathcal{L}\left(\frac{2\sinh(\frac{b}{2})}{(1+e^{-\frac{k}{2}}) e^{\frac{b}{2}}}\right)\Bigg\}=\frac{\pi^2}{2}
\end{equation}
where the sum is over the simple length spectrum $SS_T$ and $\mathcal{L}(z)$ is the Roger's dilogarithm function.
\end{theorem}

By taking the limit as $k \rightarrow 0$ in the proof of Theorem \ref{thm:oneholed}, we obtain the following identity for hyperbolic tori with one cusp.

\begin{theorem}
\label{thm:oncepuncturetorus}
For a once-punctured torus $T$ equipped with a complete hyperbolic structure of finite area,

\begin{equation}
\label{newpunctureid}
\sum_{b}\left\{\mathcal{L}\left(\operatorname{sech}^2(\frac{b}{2})\right) + 2\mathcal{L}\left(\frac{1 + e^{-b}}{2}\right)-2 \mathcal{L}\left(\frac{1 - e^{-b}}{2}\right)\right\}=\frac{\pi^2}{2}
\end{equation}
where the sum  is over the simple length spectrum $SS_T$.

\end{theorem}

Compare the above result with McShane's identity, which states that, with the same conditions,

\begin{equation}
\label{eq:McShane}
\sum_{b}\frac{1}{1+e^{b}}=\frac{1}{2}.
\end{equation}
Theorem \ref{thm:oncepuncturetorus} has an equivalent form in terms of traces and representations into
$\PSL(2,\mathbb{R})$ (see Theorem \ref{thm:puncturenew}).

\medskip
Recall that an embedded, {\em geometric} three-holed sphere/one-holed torus (collectively, simple surface) in a hyperbolic surface $S$ is an embedding from a three-holed sphere $P$/one-holed torus $T$ into $S$ such that the boundaries are mapped to geodesics.
The Luo-Tan identity for closed hyperbolic surfaces in  \cite{LuoTan11} is given by:

\begin{theorem}[Theorem 1.1 of \cite{LuoTan11}]\label{thm:LuoTan}
Let $S$ be a closed, orientable hyperbolic surface of genus $\textsf{g} \ge 2$. There exist functions $f$ and $g$ involving the Roger's dilogarithm function of the lengths of the simple closed geodesics on a three-holed sphere or a one-holed torus, such that
\begin{equation}
\label{eq:originalLuoTan}
\sum_{P} f(P) + \sum_{T} g(T) = 8 \pi^2 (\textsf{g}-1)
\end{equation}
where the first sum is over all embedded geometric three-holed spheres $P \subset S$ and the second sum extends over all  embedded geometric one-holed tori $T \subset S$.
\end{theorem}

As a corollary of Theorem \ref{thm:oneholed} we can express the Luo-Tan identity as a sum over all immersed geometric three-holed spheres $P$ of $S$ such that the restriction to $int(P)$, the interior of $P$ is an embedding into $S$.

We define a {\em quasi-embedded geometric} three-holed sphere in $S$ (orientable) to be an immersion of a three-holed sphere $P$ into $S$ which is injective on the interior $int(P)$ of $P$ such that the boundaries are mapped to geodesics, but two of the boundaries are mapped to the same geodesic. Thus, a quasi-embedded geometric three-holed sphere is contained in a unique embedded geometric one-holed torus, and conversely, every embedded geometric one-holed torus together with a non-peripheral simple closed geodesic on the torus determines a quasi-embedded geometric three-holed sphere. We have:

\begin{corollary}
\label{cor:newform}
Let $S$ be a closed, orientable hyperbolic surface of genus $\textsf{g} \ge 2$. There exist functions $f_1$ and $f_2$ involving the dilogarithms of the lengths of the simple closed geodesics in a hyperbolic three-holed sphere such that
\begin{equation}
\label{eq:newformLuoTan}
\sum_{P_1} f_1(P_1) + \sum_{P_2} f_2(P_2) = 8 \pi^2 (\textsf{g}-1)
\end{equation}
where the first sum is over all  embedded geometric three-holed spheres $P_1 \subset S$ and the second sum is over all quasi-embedded geometric three-holed spheres $P_2$ of $S$.
\end{corollary}

The explicit formulae for $f_1$ and $f_2$ are given in \S \ref{ss:pseudo}, $f_1$ is in fact just the function $f$ in Theorem \ref{thm:LuoTan}. Note that the main advantage of this reformulation is that in Theorem \ref{thm:LuoTan}, the function $g$ in \eqref{eq:originalLuoTan} was expressed in terms of an infinite sum and so the identity \eqref{eq:originalLuoTan} really involves double sums, in this formulation, the identity \eqref{eq:newformLuoTan} is expressed as a simple sum over all geometric three-holed spheres which are embeddings when restricted to the interior.

\medskip

\noindent {\bf Remark.}We can also easily get the reformulation of the Luo-Tan identity for non-closed hyperbolic surfaces whose Euler characteristics are strictly less than $-1$ as in Theorem 1.2 of \cite{LuoTan11} (Arxiv version). 

\medskip
The rest of this paper is organized as follows. In \S \ref{ss:notion}, we set the notation and apply the main theorem in \cite{LuoTan11} to a hyperbolic four-holed sphere $Q$ with four geodesic boundaries all of length $c$, to obtain an identity for $Q$. In \S \ref{ss:proofs}, we describe the relation between hyperbolic four-holed spheres $Q$ with geodesic boundaries all of length $c$, and hyperbolic one-holed tori $T$ with one geodesic boundary of length $k=2c$. We then prove Theorems \ref{thm:oneholed},  and \ref{thm:oncepuncturetorus}  by transcribing the identity for $Q$ to the one for $T$. Finally, in \S \ref{ss:pseudo}, we prove Corollary \ref{cor:newform} giving the explicit formulas for $f_1$ and $f_2$.

\bigskip

\noindent {\it Acknowledgements.} We would like to thank Martin Bridgeman, Feng Luo and Hugo Parlier for helpful discussions and comments. 

\section{Notation and results for four-holed spheres}
\label{ss:notion} We establish the Luo-Tan identity for hyperbolic four-holed spheres with boundaries of equal lengths in this section.

\noindent {\bf The Roger's Dilogarithm. }
Recall that the  Roger's Dilogarithm function $\mathcal{L}(z)$  for $z\le 1$  is defined by
  \[
  \mathcal{L}(z):=-\frac{1}{2} \int_0^z(\frac{\log u }{1-u}+\frac{\log (1-u)}{u}) \, du
  \]
see \cite{Roger1907}.   $\mathcal{L}(z)$ can also be thought of as a normalization of the regular dilogarithm function
$$Li_2(z) = \sum_{n=1}^\infty \frac{z^n}{n^2}\qquad \mbox{ for } |z| < 1,$$
where
$${\mathcal L}(z) = Li_2(z) + \frac{1}{2}\log|z|\log(1-z).$$

${\mathcal L}(z)$  is an increasing function on $z \le 1$ and it enjoys many remarkable properties, see for example \cite{Kirillov89,Leonard91,GeRiRa99} for a detailed discussion. Here we mention some of the more useful properties.

  We have $\mathcal{L}(0)=0$, $\mathcal{L}(\frac{1}{2})=\frac{\pi^2}{12}$, $ \mathcal{L}(1) =\frac{\pi^2}{6}$ and $\mathcal{L}(-1) =-\frac{\pi^2}{12}$. The Euler relations are given by:
  \[
  \mathcal{L}(x)+\mathcal{L}(1-x) = \mathcal{L}(1)=\frac{\pi^2}{6}
  \]
   for $0 \le x \le 1$ and
   \[
  \mathcal{L}(-x)+\mathcal{L}(-x^{-1}) = 2\mathcal{L}(-1)=-\frac{\pi^2}{6}
  \]
  for $x>0$. The Landen's identity is
  \[
  \mathcal{L}(\frac{-x}{1-x})=-\mathcal{L}(x)
  \]
  for $0\le x \le 1$. Note $\displaystyle \lim_{x \to -\infty}\mathcal{L}(x) = -\mathcal{L}(1) = -\frac{\pi^2}{6}$. The pentagon relation is as follows: for $x, y \in [0,1]$ and $xy \neq 1$,
      \[
       \mathcal{L}(x)+\mathcal{L}(y)+\mathcal{L}(\frac{1-x}{1-xy}) + \mathcal{L}(\frac{1-y}{1-xy}) =\mathcal{L}(xy)+\frac{\pi^2}{3}.
      \]

\medskip

\noindent {\bf The Lasso function}. The Lasso function $La(x,y)$ in \cite{LuoTan11}, which corresponds to the measure of a certain subset of the unit tangent bundle of a three-holed sphere is defined in terms of the Roger's Dilogarithm function as follows:
\begin{equation}
\label{eq:lassofunction}
La(x,y):=\mathcal{L}(y)+\mathcal{L}(\frac{1-y}{1-xy})-\mathcal{L}(\frac{1-x}{1-xy})
\end{equation}
for $0 < x < y < 1$.

\medskip

\noindent {\bf Remark.} In fact, the function $La(x,y)$ is well-defined on the square domain $0 \le x,y \le 1$. Here it is enough to focus on the upper triangle region.

\medskip
\noindent {\bf Orthogeodesics}. Given a (cone)-hyperbolic surface $S$, an {\it orthogeodesic} from $A$ to $B$ is a closed geodesic arc $\alpha$ from the set $A$ to the set $B$ which is orthogonal to $A$ and $B$. Here we usually consider $A$ and $B$ to be geodesics, however, they may be points on $S$ or cone singularities on $S$ (where the orthogonality condition there is vacuous), or they may be horocycles around cusps. The orthogeodesic is {\it simple} if the restriction to the interior is injective.

\medskip

\begin{figure}[htb]
\centering
\def\svgwidth{200pt}
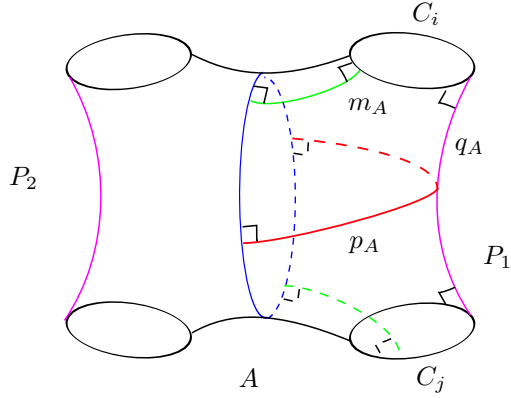
\caption{A four-holed sphere Q}
\label{fig:4sphere}
\end{figure}

\noindent{\bf The Luo-Tan identity for four-holed spheres}. Let $Q$ be a hyperbolic four-holed sphere with boundary geodesics $C_1,C_2,C_3,C_4$ all of length $c > 0$ and let $A$ be a non-peripheral simple closed geodesic on $Q$, see Figure \ref{fig:4sphere}. We obtain two isometric three-holed spheres $P_1$, $P_2$ by cutting the four-holed sphere $Q$ along $A$. Suppose that $P_1$ has boundaries $C_i$, $C_j$ and $A$, where $\{i,j\} \subset \{1,2,3,4\}$. Let $m_A$ be the length of the simple orthogeodesic from $C_i$ to $A$ in $P_1$, $p_A$ be the length of the simple orthogeodesic from $A$ to itself in $P_1$ and $q_A$ the length of the simple orthogeodesic from $C_i$ to $C_j$ in $P_1$, see Figure \ref{fig:4sphere}. Note that by  symmetry, $m_A$ is also the length of the simple orthogeodesic from $C_j$ to $A$ and that the lengths of the corresponding simple orthogeodesics on $P_2$ are given by the same quantities.

  Let $U(Q)$ be the set of unit tangent vectors over $Q$ and $\mu$ be the measure on $U(Q)$ which is invariant under the geodesic flow. We have $\mu(U(Q))=-4 \pi^2 \chi(Q) =8 \pi^2$ and applying the Luo-Tan identity to $Q$ gives:

\begin{lemma}
\label{le:4boundary}
For a hyperbolic four-holed sphere $Q$ with boundary geodesics all of length $c >0$,
\begin{equation}
\label{eq:fournew}
\sum_{A} \Bigg\{\mathcal{L}\left(\tanh^2 (\frac{p_A}{2})\right) + 2 \mathcal{L}\left(\tanh^2 (\frac{m_A}{2})\right) -
2 La\left(e^{-c}, \tanh^2(\frac{m_A}{2})\right) \Bigg\} =\frac{\pi^2}{2}
\end{equation}
where the sum is over all non-peripheral simple closed geodesics $A$ on $Q$ and
for every non-peripheral simple closed geodesic $A$, $m_A$ and $p_A$ are the lengths
of the simple orthogeodesics as described above (see Figure \ref{fig:4sphere}).
\end{lemma}

  Note that  $e^{-c} < \tanh^2(\frac{m_A}{2})$ in the above since for a fixed $c>0$, there is a unique real $m'$ such that $ \sinh(m')\sinh(\frac{c}{2}) = 1$ and for any hyperbolic structure on $Q$ and every non-peripheral simple closed geodesic $A$ on $Q$, $m' < m_A $, from which the inequality follows.

  Furthermore, we can use hyperbolic identities and the properties of the Roger's dilogarithm function $\mathcal{L}(z)$ given at the beginning of the section to obtain an equivalent form of the above identity in terms of just the quantities $\ell(A)$ and $c$. We omit the elementary but somewhat tedious details. We have the following reformulation of identity (\ref{eq:fournew}) for $Q$:
\begin{equation}
\label{eq:foursimpleid}
\sum_{A} \Bigg\{ \mathcal{L}\left(\frac{\cosh(c)+1}{\cosh(c)+ \cosh(\frac{\ell(A)}{2})}\right)+ 2\mathcal{L}\left(\frac{\cosh(\frac{c}{2}+\frac{\ell(A)}{4})}{\cosh(\frac{c}{2}) e^{\frac{\ell(A)}{4}} }\right)
-2 \mathcal{L}\left(\frac{2\sinh(\frac{\ell(A)}{4})}{(1+e^{-c})e^{\frac{\ell(A)}{4}}}\right)\Bigg\} =\frac{\pi^2}{2}
\end{equation}
where the sum extends over all non-peripheral simple closed geodesics $A$ on $Q$ and $\ell(A)$ denotes the length of $A$.

Taking the limit as $c$ tend to $0$, we obtain:
\begin{lemma}
\label{le:4cusped}
For a hyperbolic quadruply-punctured sphere $Q$,
\begin{equation}
\label{eq:ideal4sphere}
\sum_{A} \left\{\mathcal{L}\left(\operatorname{sech}^2(\frac{\ell(A)}{4})\right)+2\mathcal{L}\left(\frac{1+e^{-\frac{\ell(A)}{2}}}{2}\right) -2 \mathcal{L}\left(\frac{1-e^{-\frac{\ell(A)}{2}}}{2}\right) \right\} = \frac{\pi^2}{2}
\end{equation}
where the sum is over all simple closed geodesics $A$ on $Q$ and $\ell(A)$ denotes the length of $A$.
\end{lemma}

\section{Proofs of Theorems \ref{thm:oneholed} and \ref{thm:oncepuncturetorus}}
\label{ss:proofs}
In the following we use the widely known relationship between hyperbolic structures on one-holed tori and four-holed spheres. We claim that a hyperbolic one-holed torus with geodesic boundary of length $2c$ is paired with a hyperbolic four-holed sphere with geodesic boundaries all of length $c$, via a (two-to-one and four-to-one respectively) branched covering over a hyperbolic orbifold of genus $0$ with three cone points of angle $\pi$ and one geodesic boundary of length $c$. We provide some details below for the convenience of the reader.

\begin{figure}[htb]
\centering
\def\svgwidth{300pt}
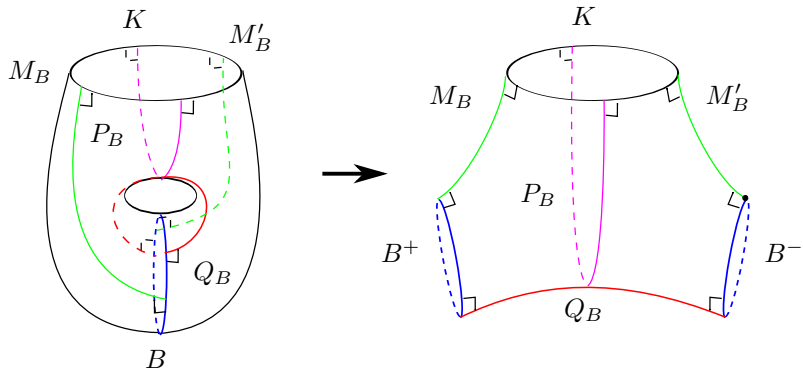
\caption{A one-holed torus $T$}
\label{fig:1holedtorus}
\end{figure}

\medskip
\noindent{\bf One-holed torus and the elliptic involution}. Let $T$ be a hyperbolic one-holed torus with boundary geodesic $K$ of length $k$ and let $B$ be a non-peripheral simple closed geodesic on $T$ with length $b$. Let $P_B$ be the simple orthogeodesic from $K$ to itself which is disjoint from $B$ and let $p_B$ be its length. Let $Q_B$ be the simple orthogeodesic from $B$ to itself which intersects $B$ at only its end points, and let $q_B$ be its length. Finally, let $M_B$ and $M_B'$ be the simple orthogeodesics from $K$ to $B$ which intersect $B$ only at the endpoint, they both have the same length which we denote by $m_B$. These geodesics and orthogeodesics on $T$, as well as on the three-holed sphere obtained by cutting $T$ along $B$ are shown in Figure \ref{fig:1holedtorus}.

Let $p_1$ and $p_2$ be the (not necessarily distinct) points on $B$ which are the endpoints of $Q_B$. Note that the two points coincide if and only if the twist parameter of the Fenchel-Nielsen coordinates with respect to the curve $B$ is an integer multiple of $b$. Generically, $p_1$ and $p_2$  separate $B$ into two arcs, we denote the midpoints of these two arcs by $\omega_1$ and $\omega_2$. Note that $\omega_1$ and $\omega_2$ are diametrically opposite on $B$.

Let $C$ be any other simple closed geodesic on $T$ which intersects $B$ exactly once. Then $C$ intersects $P_B$ exactly once at the midpoint of $P_B$ which we denote by $\omega_3$ and it intersects $B$ at either $\omega_1$ or $\omega_2$. Furthermore, these points of intersection of $C$ with $P_B$ and $B$ are diametrically opposite on $C$, and the point missed by $C$ is the midpoint of $P_C$ which is the simple orthogeodesic from $K$ to itself which is disjoint from $C$.

As every simple closed geodesic on $T$ can be obtained from a basic triple of simple closed geodesics which pair-wise intersect once, by moving along a binary tree, we see that every simple closed geodesic passes through exactly two of the three points $\omega_1, \omega_2, \omega_3$, which are diametrically opposite on the geodesic. The three points $\omega_1, \omega_2, \omega_3$ are the Weierstrass points of $T$ and rotation by $\pi$ at any one of $\omega_i$ is an isometry $\iota$ of $T$ which fixes all three points. $\iota$ is the elliptic involution of $T$, it fixes every non-peripheral simple closed geodesic, but reverses the direction.

\begin{figure}[htb]
\centering
\def\svgwidth{350pt}
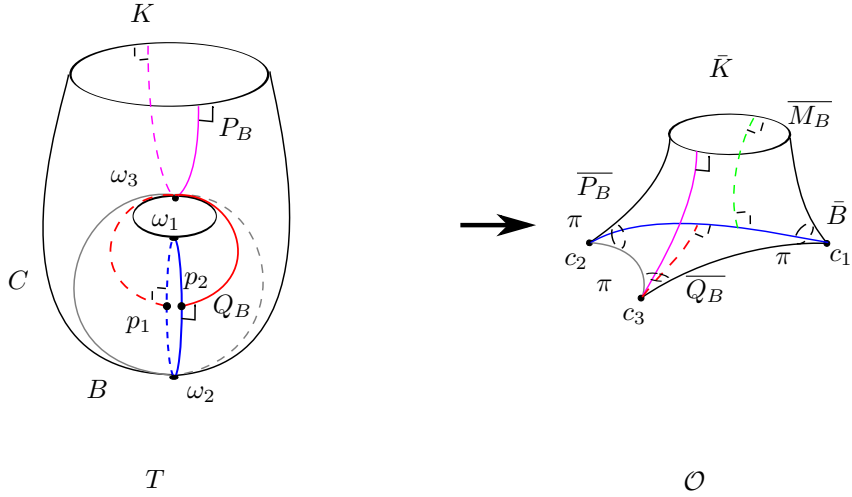
\caption{One-holed tori with orbifolds}
\label{fig:torusorbifold}
\end{figure}

We have $T/\langle \iota\rangle $ is a hyperbolic orbifold $\mathcal{O}$ of genus $0$ with three cone points $c_1, c_2, c_3$ of cone angle $\pi$ and one boundary $\overline K$  with length $\frac{k}{2}$. $T$ is a double branched cover of $\mathcal{O}$ and the Weierstrass points $\omega_1, \omega_2, \omega_3$ cover the cone points $c_1, c_2, c_3$ respectively, see Figure \ref{fig:torusorbifold}.

The boundary $K$ double covers the boundary $\overline K$ of ${\mathcal O}$,
the geodesic $B$ double covers a simple geodesic arc $\overline B$ on $\mathcal{O}$ joining $c_1$ and $c_2$, so $b=2\ell(\overline{B})$. $P_B$ double covers a simple orthogeodesic $\overline{P_B}$ from $c_3$ to the boundary disjoint from $\overline B$, so $p_B=2\ell(\overline{P_B})$. $M_B$ and $M_B'$ (as in Figure \ref{fig:1holedtorus}) cover the simple orthogeodesic $\overline{M_B}$ on ${\mathcal O}$ from $\overline K$ to $\overline B$, whose interior is disjoint from $\overline B$, so $m_B=\ell(\overline{M_B})$. $Q_B$ double covers the simple orthogeodesic $\overline{Q_B}$ from $c_3$ to $\overline B$, whose interior is disjoint from $\overline B$, so $q_B=2\ell(\overline{Q_B})$.

Reversing the construction, and keeping track of the covering data, we see every such hyperbolic orbifold ${\mathcal O}$ with boundary geodesic of length $k/2$ is double branch covered by a unique hyperbolic torus $T$ with geodesic boundary of length $k$.

\medskip

\noindent{\bf Four-holed spheres and the $\Z_2 \times \Z_2$ action}. The situation for a hyperbolic four-holed sphere $Q$ with four boundary geodesics $C_1,C_2,C_3,C_4$ of equal length $c$ is similar but slightly more complicated.

\begin{figure}[htb]
\centering
\def\svgwidth{400pt}
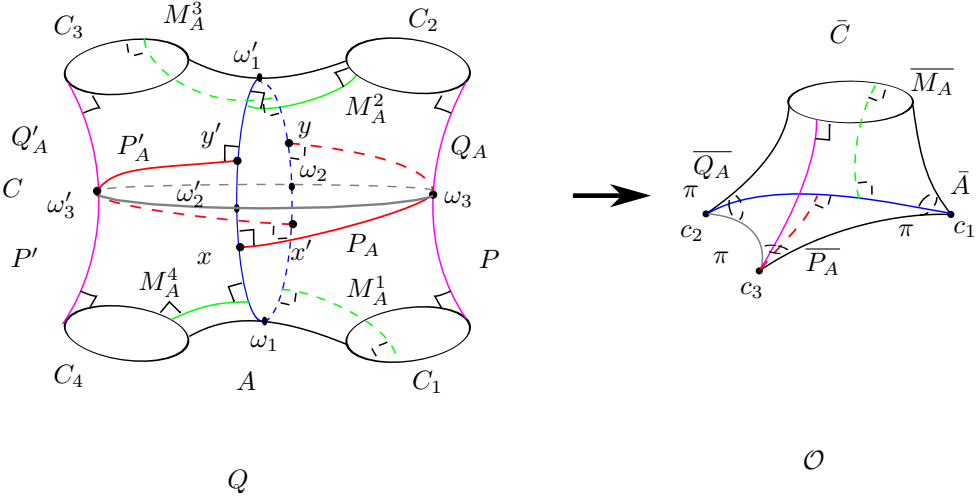
\caption{Four-holed spheres with orbifolds}
\label{fig:4holedorbifold}
\end{figure}

Let $A$ be a non-peripheral simple closed geodesic on $Q$ with length $a$. Cutting along $A$ gives two isometric three-holed spheres $P$ and $P'$, without loss of generality, assume that $C_1, C_2 \subset \partial P$ and $C_3, C_4 \subset \partial P'$. Let $P_A$ (resp. $ P_A'$) be the simple orthogeodesic on $P$ (resp. $P'$) from $A$ to itself with length $p_A$ and let $x,y$ (resp. $x',y'$) be the endpoints of $P_A$ (resp. $P_A'$) on $A$, so that $x,x',y,y'$ are in cyclic order on $A$. The midpoints of the four intervals obtained divide $A$ into four segments of equal length, denote these points by $\omega_1, \omega_2, \omega_1', \omega_2'$ respectively, where $\omega_i$ and $\omega_i'$ are diametrically opposite on $A$.. Let $M_A^i$, $i=1,2,3,4$ be the simple orthogeodesic from $C_i$ to $A$ in $P$ and $P'$, they have equal length, denoted by $m_A$ and let $Q_A$ (resp. $Q_A'$) be the simple orthogeodesic from $C_1$ to $C_2$ (resp. $C_3$ to $C_4$) in $P$ (resp. $P'$), with length $q_A$.  See Figure \ref{fig:4holedorbifold}.

Now let $C$ be any simple closed geodesic which intersects $A$ exactly twice. Then, $C$ intersects $Q_A$ and $Q_A'$ exactly once at their midpoints, which we denote by $\omega_3$, $\omega_3'$ respectively. Furthermore, $C$ intersects $A$ in either the pair of points $\omega_1, \omega_1'$ or the pair $\omega_2,\omega_2'$, and these together with $\omega_3, \omega_3'$ are equidistributed on $C$. As in the one-holed torus case, we can in fact show that every non-peripheral simple closed geodesic on  $Q$ passes through exactly two of the pairs $\{\omega_i, \omega_i\}$, which are equally distributed on the geodesic. These six points are the Weierstrass points of $Q$. A $\pi$ rotation about say $\omega_1$ is an isometry of $Q$ taking $P$ to $P'$, rotation about $\omega_1'$ induces the same involution on $Q$, which we denote by $\iota_1$. We define $\iota_2, \iota_3$ similarly, together, they generate a $\Z_2 \times \Z_2$ action on $Q$ since $\iota_1\iota_2=\iota_3$ etc., so $\langle \iota_1, \iota_2\rangle \cong Z_2 \times \Z_2$.

%

$Q /\langle \iota_1, \iota_2 \rangle$ is a genus $0$ hyperbolic orbifold ${\mathcal O}$ with three cone points $c_1, c_2, c_3$ of cone angle $\pi$ and one boundary $\overline{C}$ with length $c$. $Q$ is a quadruply branched cover of ${\mathcal O}$ with $\omega_i$ and $\omega_i'$ projecting to the same cone singularity $c_i$ ($i=1,2,3$) on ${\mathcal O}$. The four peripheral geodesics $C_i$ on $Q$ cover the peripheral geodesic $\overline{C}$ on ${\mathcal O}$ (so they have the same length), and $A$ is a quadruple cover of a simple geodesic $\overline A$ from $c_1$ to $c_2$ in ${\mathcal O}$, so $a=4 \ell({\overline A})$. $M_A^i$ $i=1,2,3,4$ cover the same simple orthogeodesic $\overline{M_A}$ on ${\mathcal O}$ from ${\overline C}$ to ${\overline A}$, so $m_A=\ell(\overline{M_A})$. $P_A$ and $P_A'$ each double covers the same simple orthogeodesic $\overline{P_A}$ on ${\mathcal O}$ from $c_3$ to $\overline A$, so $p_A=2\ell(\overline{P_A})$, and $Q_A$ and $Q_A'$ double covers the same simple orthogeodesic ${\overline Q_A}$ on ${\mathcal O}$ from $c_3$ to ${\overline C}$ which is disjoint from $\overline A$, so that $q_A=2 \ell(\overline{Q_A})$.


Again, by reversing the construction, and keeping track of the covering data, we see every hyperbolic orbifold ${\mathcal O}$ with three cone points $c_1,c_2,c_3$ of angle $\pi$ and one geodesic boundary of length $c$  is quadruply branch covered by a unique hyperbolic four-holed sphere $Q$ with geodesic boundaries all  of length $c$.

\medskip

\noindent {\bf Tori and four-holed spheres}. To summarize, every hyperbolic one-holed torus $T$ with geodesic boundary of length $2c$ is a double branched cover over a unique hyperbolic orbifold ${\mathcal O}$ with three cone points of angle $\pi$ and one geodesic boundary of length $c$, which in turn is  quadruply branched covered by a unique hyperbolic four-holed sphere $Q$ with four geodesic boundaries all of length $c$.

Furthermore, if $T$ and $Q$ are the  hyperbolic one-holed torus and four-holed sphere covering the same orbifold ${\mathcal O}$, a non-peripheral simple closed geodesic $B$ on $T$ double covers a simple geodesic arc $\gamma$ from $c_i$ to $c_j$ on $\mathcal{O}$, which in turn is quadruply covered by a unique non-peripheral simple closed geodesic $A$ on $Q$. Moreover, the simple orthogeodesic $M_B$ (or $M_B'$) on $T$ covers the simple orthogeodesic on $\mathcal{O}$ from $\overline{K}$ to $\gamma$, which in turn is covered by the simple orthogeodesic $M_A^1$ (or $M_A^2$, $M_A^3$, $M_A^4$) on $Q$; the simple orthogeodesic $Q_B$ on $T$ double covers the simple orthogeodesic on $\mathcal{O}$ from $c_l$ to $\gamma$, which in turn doubly covered by the simple orthogeodesic $P_A$ (or $P_A'$) on $Q$; and the simple orthogeodesic $P_B$ on $T$ double covers the simple orthogeodesic on $\mathcal{O}$ from $c_l$ to the boundary which is disjoint from $\gamma$, which is in turn doubly covered by the simple orthogeodesic $P_A$ (or $P_A'$) on $Q$. Therefore there is a one to one correspondence between the set of non-peripheral simple closed geodesics $\{B\}$ on $T$ and the set of non-peripheral simple closed geodesics $\{A\}$ on $Q$, and if the simple closed geodesic $B$ on $T$ corresponds to the simple closed geodesic $A$ on $Q$, we have $\ell(A)=2\ell(B)$, $m_A=m_B$, $p_A=q_B$ and $q_A=p_B$.

\medskip

\noindent {\bf Identities for one-holed tori}. From the discussion above, we see that the identity \eqref{eq:fournew} in Lemma \ref{le:4boundary} induces the following identity on the hyperbolic one-holed torus:
\begin{theorem}
\label{thm:oneholednew}
For a hyperbolic one-holed torus $T$ with boundary geodesic $K$ of length $k>0$, we have
\begin{equation}
\label{eq:1holednon0}
\sum_{B} \left\{\mathcal{L}\left(\tanh^2 (\frac{q_B}{2})\right)+2\mathcal{L}\left(\tanh^2(\frac{m_B}{2})\right) - 2 La\left(e^{-\frac{k}{2}}, \tanh^2(\frac{m_B}{2})\right) \right\}=\frac{\pi^2}{2}
\end{equation}
where $B$ is over all non-peripheral simple closed geodesics on $T$. For every non-peripheral simple closed geodesic $B$ on $T$, $q_B$ and $m_B$ are the lengths of the simple orthogeodesics shown in Figure \ref{fig:1holedtorus}.
\end{theorem}

 Note for a fixed $k>0$, there is a unique real $m''$ such that $\sinh(m'') \sinh(\frac{k}{4}) = 1$. For any hyperbolic structures on $T$ and every non-peripheral simple closed geodesics $B$ on $T$, $m'' < m_B$, thus we have $ e^{-\frac{k}{2}} < \tanh^2(\frac{m_B}{2})$.
Theorem \ref{thm:oneholednew}   in turn is equivalent to Theorem \ref{thm:oneholed} by using the reformulation  \eqref{eq:foursimpleid} of \eqref{eq:fournew}. \qed

\medskip

Let $k \to 0$, we get a once-punctured torus T. The relation between once-punctured tori and quadruply-punctured spheres applies, thus from the identity \eqref{eq:ideal4sphere} for a hyperbolic quadruply-punctured sphere, we have Theorem \ref{thm:oncepuncturetorus}. \qed

\medskip

 There is an equivalent interpretation of Theorem \ref{thm:oncepuncturetorus}, which is motivated by Bowditch \cite{Bowditch96}. Let $\Gamma$ be the fundamental group of T which a free group of rank two i.e. $\Gamma=\langle \alpha, \beta \rangle$. Any complete hyperbolic structure on a once-punctured torus T is determined by a discrete faithful and type-preserving (trace of commutator of a generating pair for $\pi_1(T)$  is $-2$)  representation $\rho:\Gamma \to \PSL(2,\mathbb{R})$. Define $g' \sim g$ in $\Gamma$ if and only if $g'$ is conjugate to $g$ or $g^{-1}$, thus $\Gamma/\!\!\sim$ represents the set of free homotopy classes of (unoriented) closed curves on T. Let $\Omega \subset \Gamma /\!\!\sim$ be the set of free homotopy classes of essential (i.e. non-trivial and non-peripheral) simple closed curves on T. For each free homotopy class $X \in \Omega$, there exists a unique non-peripheral simple closed geodesic $B$ on T representing $X$ and we have
\begin{equation}
\label{eq:basichyprep}
\tr^2(\rho(g))=4\cosh^2(\frac{\ell(B)}{2})
\end{equation}
where $g \in \Gamma$ is an arbitrary element representing the homotopy class $X$ and $\ell(B)$ denotes the length of $B$. We define a function $\phi$ on $\Omega$ by $\phi(X):=x=\tr^2(\rho(g))$ where $g \in \Gamma$ represents the homotopy class $X$. Then we have

\begin{theorem}
\label{thm:puncturenew}
For any discrete faithful, type-preserving representation $\rho:\Gamma \to \PSL(2,\mathbb{R})$,
\begin{equation}
\label{eq:newoncepuncture}
\sum_{X \in \Omega} \left\{ \mathcal{L}(\frac{4}{x}) + 2\mathcal{L}(\frac{x}{x+\sqrt{x^2-4x}})
-2\mathcal{L}(\frac{\sqrt{x^2-4x}}{x+\sqrt{x^2-4x}})\right\} = \frac{\pi^2}{2}
\end{equation}
where the sum is over all free homotopy classes in $\Omega$.
\end{theorem}

\section{Reformulation of the Luo-Tan identity}
\label{ss:pseudo}
In the following, we reformulate the function $g$ of Theorem 1.1 in \cite{LuoTan11}, which is
\begin{equation}
\label{eq:originalG}
g(T)=4 \pi^2-\sum_B 8\left[2 La\left(e^{-\ell(B)}, \tanh^2(\frac{m_B}{2})\right) +\mathcal{L}\left(\operatorname{sech}^2 (\frac{p_B}{2})\right)\right]
\end{equation}
where the sum is over all non-peripheral simple closed geodesics $B$ on $T$, and for each non-peripheral simple closed geodesic $B$, $\ell(B)$ is the length of $B$, $m_B$ and $p_B$ are the lengths of the simple orthogeodesics as in Figure \ref{fig:1holedtorus}. Here the inequality $e^{-\ell(B)}<\tanh^2(\frac{m_B}{2})$ holds for every non-peripheral simple closed geodesic $B$ on $T$. By using \eqref{eq:1holednon0} to write $4 \pi^2$ as a sum over $B \in SG_T$, we have:
\begin{align*}
g(T)=\sum_B  8\Bigg\{ & \mathcal{L}\left(\tanh^2 (\frac{q_B}{2})\right) + 2\mathcal{L}\left(\tanh^2 (\frac{m_B}{2})\right)  -\mathcal{L}\left(\operatorname{sech}^2(\frac{p_B}{2})\right)\\
& - 2 La\left(e^{-\ell(B)}, \tanh^2(\frac{m_B}{2})\right) - 2La\left(e^{-\frac{k}{2}}, \tanh^2(\frac{m_B}{2})\right)\Bigg\} \numberthis \label{eq:newgfun}
\end{align*}
where the sum is over all non-peripheral simple closed geodesics $B\in SG_T$, and for each non-peripheral simple closed geodesic $B$, $\ell(B)$ is the length of $B$, $m_B$, $p_B$ and $q_B$ are the lengths of the simple orthogeodesics $M_B, P_B$ and $Q_B$ respectively as in Figure \ref{fig:1holedtorus}.

  For any hyperbolic surface $S$ and an embedded geometric 3-holed sphere $P \subset S$, we denote the boundary geodesics by $L_1, L_2, L_3$ with  lengths  $l_1,l_2,l_3$, respectively. For $\{i, j, k\} = \{1,2,3\}$, let $M_i$ be the shortest geodesic arc between $L_j$ and $L_k$ with its length $m_i$ and $B_i$ be the shortest non-trivial geodesics arc from $L_i$ to itself with its length $p_i$. We define
\begin{align*}
f_1(P)&:=4 \pi^2 -  8\left[\sum_{i=1}^3 \left(\mathcal{L}\left(\operatorname{sech}^2(\frac{m_i}{2})\right) + \mathcal{L}\left(\operatorname{sech}^2(\frac{p_i}{2})\right)\right) + \sum_{i \neq j} La\left(e^{-l_i},\tanh^2(\frac{m_j}{2})\right) \right]\\
&:= 8\left[\sum_{i=1}^3 \left(\mathcal{L}\left(\tanh^2(\frac{m_i}{2})\right) - \mathcal{L}\left(\operatorname{sech}^2(\frac{p_i}{2})\right)\right)  -\sum_{i \neq j} La\left(e^{-l_i},\tanh^2(\frac{m_j}{2})\right) \right].  \numberthis \label{eq:newf1fun}
\end{align*}

 For a quasi-embedded geometric 3-holed sphere $P$ of $S$, let $B \subset S$ be the simple closed geodesic  which is the image of two of the boundaries of $P$, and $T$ the one-holed torus which is the image of $P$ in $S$. We denote  the length of $B$ by $\ell(B)$, and let $m_B$, $p_B$ and $q_B$ the lengths of the simple orthogeodesics $M_B, P_B$ and $Q_B$  as in Figure \ref{fig:1holedtorus}. We define
\begin{align*}
f_2(P):=8\Bigg[&\mathcal{L} \left(\tanh^2 (\frac{q_B}{2})\right) + 2\mathcal{L}\left(\tanh^2 (\frac{m_B}{2})\right) - \mathcal{L}\left(\operatorname{sech}^2(\frac{p_B}{2})\right)\\
&- 2 La\left(e^{-\ell(B)}, \tanh^2(\frac{m_B}{2})\right) - 2La\left(e^{-\frac{k}{2}},\tanh^2(\frac{m_B}{2})\right)\Bigg]. \numberthis \label{eq:newf2fun}
\end{align*}

 Thus \eqref{eq:newgfun} becomes
\[
g(T)= \sum_{P} f_2(P)
\]
where the sum is over all quasi-properly immersed geometric 3-holed spheres $P$ of $T$. Therefore we get Corollary \ref{cor:newform} simply by reformulating the function $g$ in \eqref{eq:originalLuoTan}.

\end{document}